\documentclass[11pt]{article}
\usepackage{amsmath}
\usepackage{amssymb}
\makeatletter \@addtoreset{figure}{section} \makeatother
\makeatletter
\long\def\@makecaption#1#2{%
   \vskip 10\p@
   \setbox\@tempboxa\hbox{{#1}\ \ #2}%
   \ifdim \wd\@tempboxa >\hsize
       {#1}\ \ #2\par
   \else
       \hbox to\hsize{\hfil\box\@tempboxa\hfil}%
   \fi}
\makeatother \makeatletter
\newcommand
\myroman[1]{\@roman #1}
\makeatother

\newtheorem{thm}{Theorem}[section]
\newtheorem{cor}[thm]{Corollary}
\newtheorem{lem}[thm]{Lemma}
\newtheorem{rem}[thm]{Remark}

\newcommand{\qed}{{\hfill\rule{3pt}{7pt}}}
\def\pf{\noindent {\it Proof} }

\def\qed{\hfill \rule{4pt}{7pt}}
\def\pf{\noindent {\it Proof} }
\title{\bf Soliton solutions for  coupled Schr\"{o}dinger systems with sign-changing potential}
\author
{
Chungen Liu
\footnote
{School of Mathematical Sciences and LPMC,  Nankai University,  Tianjin 300071,  China.\newline\indent \ \
Partially supported by NNSF of China(11071127, 10621101 )
and National Basic Research Program of China(973 Program: 2011CB808002).\newline\indent \ \ E-mail: liucg@nankai.edu.cn
}~~~~~~~
~~~~~~~~~~
Youquan Zheng
\footnote
{School of Mathematical Sciences,  Nankai University,  Tianjin 300071,  China.\newline
 \indent \ \ E-mail: zhengyq@mail.nankai.edu.cn
}
}
\begin{document}
\date{}
\maketitle {\bf \noindent\large Abstract}
In this paper,
a class of coupled systems of nonlinear Schr\"{o}dinger equations with sign-changing potential,
including the linearly coupled case,
is considered.
The existence of non-trivial bound state solutions via linking methods for cones in Banach spaces is proved.

\section{Introduction and main results}
Recently,
many mathematicians focused their attention to  coupled nonlinear Schr\"{o}dinger systems.
From the viewpoint of physics,
coupled Schr\"{o}dinger systems  arise from the models of a lot of natural phenomena.
A typical example is the study of the dynamics of coupled Bose-Einstein condensates and the following equation is derived
\begin{equation}\label{e1.1}
\left
\{
\begin{array}{lll}
i\frac{\partial \psi_{1}}{\partial t}= (-\partial^{2}/\partial x^2 + V_{1} +U_{11}|\psi_{1}|^{2}+U_{12}|\psi_{2}|^2)\psi_{1} + \lambda \psi_{2},\\
i\frac{\partial \psi_{2}}{\partial t}= (-\partial^{2}/\partial x^2 + V_{2} + U_{22}|\psi_{2}|^{2} + U_{21}|\psi_{1}|^2)\psi_{2} + \lambda \psi_{1}.
\end{array}
\right.
\end{equation}
Such systems of equations also appear in nonlinear optical models and many other physical contexts,
see \cite{BGT} for detail discussions.
For such coupled systems,
the solutions of the form
$\psi_{j} = u_{j}\exp(i\omega_{j}t)$
(standing waves) are interesting,
where
$u_{j}$
solve the following system
\begin{equation}\label{e1.2}
\left
\{
\begin{array}{lll}
-\frac{\partial^{2}u_{1}}{\partial x^2} + (V_{1}+\omega_{1})u_{1} = -(U_{11}|u_{1}|^{2} + U_{12}|u_{2}|^2)u_{1} - \lambda u_{2},\\
-\frac{\partial^{2}u_{2}}{\partial x^2} + (V_{2}+\omega_{2})u_{2} = -(U_{22}|u_{2}|^{2} + U_{21}|u_{1}|^2)u_{2} - \lambda u_{1}.
\end{array}
\right.
\end{equation}

In this paper,
we will consider the following  coupled system of nonlinear Schr\"{o}dinger equations
\begin{equation}\label{e1.8}
\left
\{
\begin{array}{lll}
- \Delta u_{1} + (b_{1}(x) - \lambda V_{1}(x))u_{1} = W_{t}(x, u_{1}, u_{2})+ \lambda\gamma(x)u_{2},\\
- \Delta u_{2} + (b_{2}(x) - \lambda V_{2}(x))u_{2} = W_{s}(x, u_{1}, u_{2}) + \lambda\gamma(x)u_{1},\\
u_{1}, u_{2}\in H^{1}(\mathbf{R}^{N}),
\end{array}
\right.
\end{equation}
here and in the sequel,
$V_{i}\in L^{\infty}(\mathbf{R}^{N})$,
$\gamma\in L^{\infty}(\mathbf{R}^{N})$,
$i = 1,2$,
$\nabla_z W=(W_{t},W_{s})$
is the gradient of
$W(x, t, s)$
with respect to
$z = (t,s)\in \mathbf{R}^{2}$
and we will write
$W(x,z) = W(x,t,s)$ for convenience.
We divide our discussions into two cases.\\
{\it The non-radially  symmetric case.}
We assume
$b_{i}(x)$
satisfying the following conditions
\begin{enumerate}
\item[{(B)}]
for
$i = 1,2$,
$b_{i}\in C(\mathbf{R}^{N})$,
there exists a constant $b_{i}^{0}>0$ such that
$\displaystyle\inf_{x\in\mathbf{R}^{N}}b_{i}(x) \geq b_{i}^{0}$,
and the $n$ dimensional Lebesgue measure
$meas\{x\in \mathbf{R}^{N}|\, b_{i}(x) \leq M\} < \infty$ for any $M>0$.
\end{enumerate}
We assume
$W$
satisfying the following conditions.
\begin{enumerate}
\item[{(W$_1$)}]
$W\in C^{1}(\mathbf{R}^{N}\times \mathbf{R}^2)$,
there exists
$p\in (2, 2^{*})$
such that
$0\le W(x, z) \leq C(1 + |z|^{p})$, $\forall \,(x,z)\in \mathbf{R}^{N}\times \mathbf{R}^2$,
here,
$2^*=\frac{2N}{N-2}$
if
$N > 2$
and
$2^*=+\infty$
if
$N = 1, 2$,
\item[{(W$_2$)}]
$\displaystyle\lim_{|z|\rightarrow \infty}\frac{W(x, z)}{|z|^2} = + \infty$
uniformly for
$x\in \mathbf{R}^{N}$,
\item[{(W$_3$)}]
$W_t(x, 0, s)= 0$,
$W_s(x, t, 0) =0$
for any
$x\in \mathbf{R}^{N}$,
$s\in \mathbf{R}$, $t\in \mathbf{R}$,
and
$\displaystyle\lim_{|z|\rightarrow 0}\frac{W(x, z)}{|z|^2} = 0$
uniformly for
$x\in \mathbf{R}^{N}$,
\item[{(W$_4$)}]
set
$\mathcal{W}(x, z) = \nabla_zW(x, z)\cdot z - 2W(x, z)$,
then there exists
$\theta  \geq 1$ such that $\theta\mathcal{W}(x, z) \geq \mathcal{W}(x, \eta z)$,
$\forall\,(x, z)\in \mathbf{R}^{N} \times \mathbf{R}^2$
and
$\eta\in [0, 1]$.
\end{enumerate}
\noindent{\bf Remark.}
(1) From
(W$_4$)
and
$\mathcal{W}(x, 0)=0$,
we see that
$\mathcal{W}(x, z)\ge 0$
for any
$(x,z)\in \mathbf{R}^{N}\times \mathbf{R}^{2}$
by taking
$\eta=0$.
So we have
$\nabla_zW(x, z)\cdot z \ge 2W(x, z)$.

\noindent(2) From condition
(W$_3$),
when
$\lambda \gamma(x)\ne 0$,
$\forall x\in \mathbf{R}^{N}$,
for a non-trivial solution
${\bf u}=(u_1,u_2)$
of the problem
$(\ref{e1.8})$,
it is easy to see that
$u_1\ne 0$
and
$u_2\ne 0$,
so
${\bf u}$
does not have an immediate counterpart for a single equation.
We also remind that under the above conditions the potential
$b_i(x)-\lambda V_i(x)$
may change sign since
$\lambda\in \mathbf{R}$,
see Theorem \ref{t1.1} below.

In this case, we have the following main result.
\begin{thm}\label{t1.1}
If {\rm (B)} and {\rm (W$_1$)--(W$_4$)} hold,
the problem
$(\ref{e1.8})$
possesses a non-trivial solution for every
$\lambda\in \mathbf{R}$.
\end{thm}
{\it The radially  symmetric case.}
We assume that
$b_{i}(x)$
satisfy the following condition
\begin{enumerate}
\item[{(B)$_{r}$}]
for $i = 1,2$,
$b_{i}\in C(\mathbf{R}^{N})$,
there exists a  constant
$b_{i}^{0}>0$
such that
$\displaystyle\inf_{x\in\mathbf{R}^{N}}b_{i}(x) \geq b_{i}^{0}$,
and
$b_{i}$
are radially symmetric,
i.e.,
$b_{i}(x) = b_{i}(|x|), \forall\, x\in \mathbf{R}^{N}$,
\end{enumerate}
and
$V_{i}(x)$,
$\gamma(x)$,
$W(x,z)$
further satisfy
\begin{enumerate}
\item[{(V)$_{r}$}]
for $i = 1,2$,
$V_{i}(x) = V_{i}(|x|)$, $\gamma(x)=\gamma(|x|), \forall\, x\in \mathbf{R}^{N}$.
\end{enumerate}
\begin{enumerate}
\item[{(W$_5$)}] $W(x,z)=W(|x|,z), \,\forall\,(x,z)\in \mathbf{R}^{N}\times \mathbf{R}^2$.
\end{enumerate}
For this case we have the following result.
\begin{thm}\label{t1.2}
If {\rm (B)$_{r}$, (V)$_{r}$} and {\rm (W$_{1}$)--(W$_{5}$)} hold,
the problem
$(\ref{e1.8})$
possesses a non-trivial radially symmetric solution for every
$\lambda\in \mathbf{R}$.
\end{thm}

Next,
we consider some special cases of (\ref{e1.8}).
Firstly,
we consider some linearly coupled systems.
Precisely,
we assume that
$W_t(x, t, s)$
dose not depend on
$s$
and
$W_s(x, t, s)$
does not depend on
$t$,
that is to say one can write (\ref{e1.8}) as
\begin{equation}\label{e1.9}
\left
\{
\begin{array}{lll}
- \Delta u_{1} + (b_{1}(x) - \lambda V_{1}(x))u_{1} = f(x,u_{1}) + \lambda\gamma(x)u_{2},\\
- \Delta u_{2} + (b_{2}(x) - \lambda V_{2}(x))u_{2} = g(x,u_{2}) + \lambda\gamma(x)u_{1},\\
u_{1}, u_{2}\in H^{1}(\mathbf{R}^{N}).
\end{array}
\right.
\end{equation}
In this case,
we assume that
$f$,
$g \in C(\mathbf{R}^{N}\times \mathbf{R})$
satisfy
\begin{enumerate}
\item[{(f$_{1}$)}]
$\exists\, p_{1}\in (2, 2^{*})$
such that
$|f(x, t)| \leq C(1 + |t|^{p_{1}-1})$,
$f(x, t)t \geq 0$,
$\forall\,(x,t)\in\mathbf{R}^N\times \mathbf{R}$,
\item[{(f$_2$)}]
set
$F(x, t)=\int_{0}^{t}f(x, t){\rm d}t$,
$\displaystyle\lim_{|t|\rightarrow \infty}\frac{F(x, t)}{|t|^{2}} = + \infty$
uniformly in
$x\in \mathbf{R}^{N}$,
\item[{(f$_3$)}]
$\displaystyle\lim_{t\rightarrow 0}\frac{f(x, t)}{t} = 0$
uniformly in
$x\in \mathbf{R}^{N}$,
\item[{(f$_4$)}]
$\mathcal{F}(x, t) = f(x, t)t - 2F(x, t)$,
then there exists
$\theta_{1} \geq 1$
such that
$\theta_{1} \mathcal{F}(x, t) \geq \mathcal{F}(x, \eta t)$,
$\forall\,(x, t)\in \mathbf{R}^{N} \times \mathbf{R}$ and $\eta\in [0, 1]$,
\item[{(g$_{1}$)}]
$\exists\,p_2\in (2, 2^{*})$
such that
$|g(x, s)| \leq C(1 + |s|^{p_{2}-1})$,
$g(x, s)s \geq 0$,
$\forall\,(x,s)\in\mathbf{R}^N\times \mathbf{R}$,
\item[{(g$_2$)}]
set
$G(x, s) = \int_{0}^{s}g(x, s){\rm d}s$,
$\displaystyle\lim_{|s|\rightarrow \infty}\frac{G(x, s)}{|s|^{2}} = + \infty$
uniformly in
$x\in \mathbf{R}^{N}$,
\item[{(g$_3$)}]
$\displaystyle\lim_{s\rightarrow 0}\frac{g(x, s)}{s} = 0$
uniformly in
$x\in \mathbf{R}^{N}$,
\item[{(g$_4$)}]
$\mathcal{G}(x, s) = g(x, s)s - 2G(x, s)$,
then there exists
$\theta_{2} \geq 1$
such that
$\theta_{2} \mathcal{G}(x, s) \geq \mathcal{G}(x, \eta s)$,
$\forall\,(x, s)\in \mathbf{R}^{N} \times \mathbf{R}$ and $\eta\in [0, 1]$.
\end{enumerate}
\begin{thm}\label{t1.3}
If
{\rm (B), ({\rm f$_{1}$)--(f$_{4}$})}
and
{\rm ({\rm g$_{1}$)--(g$_{4}$})} hold,
the problem
$(\ref{e1.9})$
possesses a non-trivial solution for every
$\lambda\in \mathbf{R}$.
\end{thm}
\pf.
Set
$W(x, t, s) = F(x, t) + G(x, s)$,
it is easy to see that
(W$_{1}$)
and
(W$_{4}$)
hold.

As for
(W$_{2}$),
from
(f$_{2}$)
and
(g$_{2}$),
$\forall M > 0$,
there exists
$R>0$
such that
$\frac{F(x, t)}{|t|^2} > 2M$
when
$|t| \geq R$
and
$\frac{G(x, s)}{|s|^2} > 2M$
when
$|s| \geq R$.
Then
$$\frac{F(x, t) + G(x, s)}{t^2 + s^2} \geq \frac{F(x, t) + G(x, s)}{2\max(|t|^2, |s|^2)} > M$$
when
$\max(|t|, |s|) \geq  R$.
So
$\displaystyle\lim_{|z|\rightarrow \infty}\frac{W(x, z)}{|z|^2} = + \infty$
uniformly for
$x\in \mathbf{R}^{N}$.

From
(f$_{3}$),
(g$_{3}$)
and the continuity of
$f$
and
$g$,
we can see
$f(x, 0) = 0 = g(x, 0)$,
so
$W_t(x, 0, s)= 0$,
$W_s(x, t, 0) =0$
for any
$x\in \mathbf{R}^{N}$,
$s\in \mathbf{R}$,
$t\in \mathbf{R}$.
Also from
(f$_{3}$)
and
(g$_{3}$),
we have
$\displaystyle\lim_{|t|\to 0}\frac{F(x, t)}{|t|^2} = 0$
and
$\displaystyle\lim_{|s|\to 0}\frac{G(x, s)}{|s|^2} = 0$,
so
$$0 \leq \frac{F(x, t) + G(x, s)}{|t|^2 + |s|^2} \leq \frac{F(x, t)}{|t|^2} + \frac{G(x, s)}{|s|^2} \to 0.$$
So
(W$_{3}$)
holds.
From Theorem \ref{t1.1},
we get the assertion.
\qed

As in Theorem \ref{t1.2},
assuming that
$f(x,t)$
and
$g(x,s)$
further satisfy
\begin{enumerate}
\item[{(f$_{5}$)}]
$f(x,t) = f(|x|,t)$,
for any
$(x,t)\in \mathbf{R}^{N}\times \mathbf{R}$,
\item[{(g$_{5}$)}]
$g(x,s) = g(|x|,s)$,
for any
$(x,s)\in \mathbf{R}^{N}\times \mathbf{R}$,
\end{enumerate}
also setting
$W(x, t, s) = F(x, t) + G(x, s)$
and by the same reason as in the proof of Theorem \ref{t1.3},
we have the following consequence.
\begin{thm}\label{t1.4}
If
{\rm (B)$_{r}$, (V)$_{r}$, ({\rm f$_{1}$)--(f$_{5}$})}
and
{\it ({\rm g$_{1}$)--(g$_{5}$})}
hold,
the problem
$(\ref{e1.9})$
possesses a non-trivial radially symmetric solution for every
$\lambda\in \mathbf{R}$.
\end{thm}

By taking
$f(x,t) =c_1(x) |t|^{p_{1} - 2}t$
and
$g(x, s) = c_2(x)|s|^{p_{2} - 2}s$
with
$c_{i}\in L^{\infty}(\mathbf{R}^{N})$
and
$\displaystyle\inf_{x\in \mathbf{R}^{N}} c_i(x)>0$,
$i=1,2$,
we get the following system
\begin{equation}\label{e1.10}
\left
\{
\begin{array}{lll}
- \Delta u_{1} + (b_{1}(x) - \lambda V_{1}(x))u_{1} = c_{1}(x)|u_{1}|^{p_{1}-2}u_{1} + \lambda\gamma(x) u_{2},\\
- \Delta u_{2} + (b_{2}(x) - \lambda V_{2}(x))u_{2} = c_{2}(x)|u_{2}|^{p_{2}-2}u_{2} + \lambda\gamma(x) u_{1},\\
u, v\in H^{1}(\mathbf{R}^{N}),
\end{array}
\right.
\end{equation}
then for
$p_{1}, p_{2}\in (2,\, 2^*)$,
we have the following consequences.
\begin{cor}\label{c1.5}
If
{\it (B)}
holds,
the problem
$(\ref{e1.10})$
possesses a non-trivial solution for every
$\lambda\in \mathbf{R}$.
\end{cor}
\begin{cor}\label{c1.6}
If
{\rm (B)$_{r}$, (V)$_{r}$}
hold,
and
$c_{i}(x) = c_{i}(|x|)$
for any
$x\in \mathbf{R}^{N}$,
$i=1,2$,
the problem
$(\ref{e1.10})$
possesses a non-trivial radially symmetric solution for every
$\lambda\in \mathbf{R}$.
\end{cor}

Secondly,
by taking
$W(x, t, s) = \frac{1}{4}t^{4} + \frac{1}{2}t^2s^2 + \frac{1}{4}s^4$,
we get the following systems
\begin{equation}\label{e1.11}
\left
\{
\begin{array}{lll}
- \Delta u_{1} + (b_{1}(x) - \lambda V_{1}(x))u_{1} = u_{1}^3 + u_{2}^2u_{1} + \lambda\gamma(x) u_{2},\\
- \Delta u_{2} + (b_{2}(x) - \lambda V_{2}(x))u_{2} = u_{2}^3 + u_{1}^2u_{2} + \lambda\gamma(x) u_{1},\\
u_{1}, u_{2}\in H^{1}(\mathbf{R}^{N}).
\end{array}
\right.
\end{equation}
as consequences of Theorem \ref{t1.1} and \ref{t1.2},
we have
\begin{cor}\label{t1.7}
If
{\rm (B)}
holds,
the problem
$(\ref{e1.11})$
possesses a non-trivial solution for every
$\lambda\in \mathbf{R}$.
\end{cor}
\begin{cor}\label{t1.8}
If
{\rm (B)$_{r}$}
and
{\rm (V)$_{r}$}
hold,
the problem
$(\ref{e1.11})$
possesses a non-trivial radially symmetric  solution for every
$\lambda\in \mathbf{R}$.
\end{cor}

The study of linearly coupled Schr\"{o}dinger systems from the mathematical point of view began very recently,
see
\cite{A08,ACR07,ACR08,BGT}.
In
\cite{ACR07},
the authors proved the existence of positive ground state solution of
the following system of nonlinear Schr\"{o}dinger equations for
$0 < \lambda < 1$,
\begin{equation}\label{e1.3}
\left
\{
\begin{array}{lll}
- \Delta u + u = (1 + a(x))|u|^{p-2}u + \lambda v,\\
- \Delta v + v = (1 + b(x))|v|^{p-2}v + \lambda u,\\
u, v\in H^{1}(\mathbf{R}^{N}),
\end{array}
\right.
\end{equation}
with
$a, b \in L^{\infty}(\mathbf{R}^{N})$,
$\displaystyle\lim_{|x|\to \infty}a(x) = \displaystyle\lim_{|x|\to \infty}b(x) = 0$, $\displaystyle\inf_{\mathbf{R}^{N}}\{1 + a(x)\} > 0$,
$\displaystyle\inf_{\mathbf{R}^{N}}\{1 + b(x)\} > 0$
and
$a(x) + b(x) \geq 0$.
In \cite{ACR08},
the authors devoted to the study the multi-bump solitons of the following system
\begin{equation}\label{e1.4}
\left
\{
\begin{array}{lll}
- \Delta u + u - u^3 = \epsilon v,\\
- \Delta v + v - v^3 = \epsilon u,\\
u, v\in H^{1}(\mathbf{R}^{N}),
\end{array}
\right.
\end{equation}
in
$\mathbf{R}^{N}$
with dimension
$N = 1, 2, 3$.
In \cite{A08},
A. Ambrosetti studied the following two systems
\begin{equation}\label{e1.5}
\left
\{
\begin{array}{lll}
- u_{1}^{''} + u_{1} = (1 + \varepsilon a_{1}(x))u_{1}^{3} + \gamma u_{2},\\
- u_{2}^{''} + u_{2} = (1 + \varepsilon a_{2}(x))u_{2}^{3} + \gamma u_{1},\\
u_{1}, u_{2}\in H^{1}(\mathbf{R}),
\end{array}
\right.
\end{equation}
\begin{equation}\label{e1.6}
\left
\{
\begin{array}{lll}
-\varepsilon^2u_{1}^{''} + u_{1} + U_{1}(x)u_{1} = u_{1}^{3} + \gamma u_{2},\\
-\varepsilon^2u_{2}^{''} + u_{2} + U_{2}(x)u_{2} = u_{2}^{3} + \gamma u_{1},\\
u_{1}, u_{2}\in H^{1}(\mathbf{R}),
\end{array}
\right.
\end{equation}
and proved the existence of non-trivial solution for (\ref{e1.5}) under the conditions
$a_{i} \in L^{\infty}(\mathbf{R})$,
$\displaystyle\lim_{|x|\to \infty}a_{i}(x) = 0$,
$i = 1, 2$,
$0 < \gamma < 1, \gamma \neq 3/5$,
and (\ref{e1.6}) possesses a solution concentrating at nondegenerate stationary points of the sum
$U_{1} + U_{2}$
when
$\varepsilon \to 0$
under the conditions
$U_i\in L^{\infty}$
and
$\displaystyle\inf_{x\in \mathbf{R}}U_i(x)>-1,\;i=1,2$.
The main tools in \cite{A08,ACR07,ACR08} are the perturbation techniques,
we refer \cite{AM} for ditailed discussions about these methods.
In \cite{BGT},
the following system was considered
\begin{equation}\label{e1.7}
\left
\{
\begin{array}{lll}
- u_{1}^{''} + a(x)u_{1} - b(x)u_{2} = c(x)H_1(u_{1}, u_{2})u_{1},\\
- u_{2}^{''} + d(x)u_{2} - e(x)u_{1} = f(x)H_2(u_{1}, u_{2})u_{2},\\
u_{1}, u_{2}\in H^{1}(\mathbf{R}),
\end{array}
\right.
\end{equation}
the authors got a non-trivial solution via Krasnoselskii fixed point theory.
We note that the potentials in systems (\ref{e1.3})-(\ref{e1.7}) are positive.

To prove the main theorem,
we deal with the existence problem of non-trivial solutions by variational methods.
We first study an eigenvalue problem,
whose eigenfunctions are solutions of (\ref{e1.8}) but without the nonlinear term,
then the non-zero critical point of the functional related to
the nonlinear perturbation of this eigenvalue problem is a weak solution of (\ref{e1.8}).
To find the critical point,
we use a critical point theorem developed
by Degiovanni and Lancelotti in \cite{DL07}.

The rest of the paper is organized as follows.
The variational setting is contained in section 2.
In section 3, we study the eigenvalue problem.
We prove that there exists a divergent sequence of eigenvalues
which are defined by the cohomological index.
We prove Theorem \ref{t1.1} and \ref{t1.2} in section 4.

\section{Variational setting}
Let
$H_{1}:= \{u_1\in H^{1}(\mathbf{R}^{N})|\int_{\mathbf{R}^{N}}b_{1}(x)u_1^2{\rm d}x < \infty\}$,
then
$H_{1}$
is a Hilbert Space with inner product
$\langle u_1, v_1\rangle_{1} = \int_{\mathbf{R}^{N}}(\nabla u_1\cdot \nabla v_1 + b_{1}(x)u_1v_1){\rm d}x$
and norm
$\|u_1\|_{1}^2 = \langle u_1,u_1 \rangle_{1}$.
Similarly,
let $H_{2}:= \{u_2\in H^{1}(\mathbf{R}^{N})|\int_{\mathbf{R}^{N}}b_{2}(x)u_2^2{\rm d}x < \infty\}$,
then
$H_{2}$
is a Hilbert Space with inner product
$\langle u_2, v_2\rangle_{2} = \int_{\mathbf{R}^{N}}(\nabla u_2\cdot \nabla v_2 + b_{2}(x)u_2v_2){\rm d}x$
and norm
$\|u_2\|_{2}^2 = \langle u_2, u_2\rangle_{2}$.

For the non-radially symmetric case,
by the condition (B),
$H_{1}$
and
$H_{2}$
can be compactly embedded into
$L^{p}(\mathbf{R}^{N})$,
$2\leq p < 2^*$ (see for example, \cite{BW, W}).
Set
$\mathcal{H}:= H_{1} \times H_{2}$,
then
$\mathcal{H}$
is a Hilbert space with inner product
$\langle\cdot,\cdot\rangle = \langle\cdot, \cdot\rangle_{1} + \langle\cdot, \cdot\rangle_{2}$
and with norm
$\|{\bf u}\|^{2} = \|u_{1}\|_{1}^{2} + \|u_{2}\|_{2}^{2}$
for
${\bf u} = (u_{1}, u_{2})$.

For the radially symmetric case,
let
$H_{1,r}: = \{u_1\in H_{1}| u_1$
is radially symmetric\},
$H_{2,r}:= \{u_2\in H_{2}|u_2$
is radially symmetric\},
then
$H_{i,r}$
is a Hilbert Space with inner product
$\langle\cdot, \cdot\rangle_{i}$
and norm
$\|\cdot\|_{i}$
for
$i = 1, 2$.
By condition
(B)$_r$,
$H_{i,r}$
can be  compactly embedded into
$L^{p}(\mathbf{R}^{N})$,
$2\leq p < 2^*$
for
$i=1,2$
(see \cite{BW, W}).
In this case,
we set
$\mathcal{H}_r:= H_{1,r} \times H_{2,r}$,
then
$\mathcal{H}_{r}$
is a Hilbert space with inner product
$\langle\cdot,\cdot\rangle = \langle\cdot, \cdot\rangle_{1} + \langle\cdot, \cdot\rangle_{2}$
and with norm
$\|{\bf u}\|^{2} = \|u_{1}\|_{1}^{2} + \|u_{2}\|_{2}^{2}$
for
${\bf u} = (u_{1}, u_{2})$.

In order to prove Theorem \ref{t1.1}, we define a functional
$\Psi: \mathcal{H}\to \mathbf{R}$
by
\begin{equation}\label{e2.1}
\Psi({\bf u})= E({\bf u}) - \lambda J({\bf u}) - P({\bf u}),\;\;{\bf u} = (u_{1}, u_{2})\in \mathcal{H},
\end{equation}
where
\begin{equation}\label{e2.2}
E({\bf u})=\frac{1}{2}\|{\bf u}\|^2,
\end{equation}
\begin{equation}\label{e2.3}
J({\bf u})=\int_{\mathbf{R}^{N}}\left(\frac{1}{2}V_{1}(x)u_{1}^2 + \gamma(x) u_{1} u_{2} + \frac{1}{2}V_{2}(x)u_{2}^2\right){\rm d}x,
\end{equation}
and
\begin{equation}\label{e2.4}
P({\bf u})= \int_{\mathbf{R}^{N}}W(x, {\bf u}){\rm d}x = \int_{\mathbf{R}^{N}}W(x, u_{1}, u_{2}){\rm d}x,
\end{equation}
then these four functionals are
$C^{1}$,
and for
${\bf u}=(u_1,u_2)$,
${\bf v}=(v_1,v_2)\in \mathcal{H}$,
there hold
\begin{equation}\label{e2.5}
\langle E'({\bf u}),  {\bf v}\rangle = \int_{\mathbf{R}^{N}}\left(\nabla u_{1} \cdot \nabla v_{1}
+ b_{1}(x)u_{1}v_{1}\right){\rm d}x + \int_{\mathbf{R}^{N}}\left(\nabla u_{2} \cdot \nabla v_{2} + b_{2}(x)u_{2}v_{2}\right){\rm d}x,
\end{equation}
\begin{equation}\label{e2.6}
\langle J'({\bf u}),  {\bf v}\rangle = \int_{\mathbf{R}^{N}}\left(V_{1}(x)u_{1}v_{1} + \gamma(x) u_{2}v_{1} + \gamma(x) u_{1}v_{2}+ V_{2}(x)u_{2}v_{2}\right){\rm d}x,
\end{equation}
\begin{equation}\label{e2.7}
\langle P'({\bf u}), {\bf v}\rangle = \int_{\mathbf{R}^{N}}\left(W_{t}(x, u_{1}, u_{2})v_{1} + W_{s}(x, u_{1}, u_{2})v_{2}\right){\rm d}x,
\end{equation}
\begin{equation}\label{e2.8}
\langle \Psi'({\bf u}), {\bf v}\rangle = \langle E'({\bf u}),  {\bf v}\rangle - \lambda \langle J'({\bf u}),  {\bf v}\rangle - \langle P'({\bf u}), {\bf v}\rangle.
\end{equation}
It is clear that critical points of
$\Psi$
are weak
solutions of
(\ref{e1.8}).

For the radially symmetric case,
we can also define these four functionals and
(\ref{e2.5})-(\ref{e2.8})
hold,
the only difference is the domain
$\mathcal{H}$
of the functional
$\Psi$
is replaced by
$\mathcal{H}_{r}$.
And the critical points of the functional
$\Psi$
are radially symmetric weak solutions of
(\ref{e1.8}).

In order to find a critical point of $\Psi$,
we need the following critical point theorem.
It was proved in \cite{DL07},
where the functional was supposed to satisfy the $(PS)$ condition.
Recently,
in \cite{D09},
the author extended it to more general case
(the functional space is completely regular topological space or metric space). As observed in \cite{LZ},
if the functional space is a real Banach space,
according to the proof of Theorem 6.10 in \cite{D09},
the Cerami condition is sufficient for the compactness of the set of critical points at a fixed level
and the first deformation lemma to hold (see \cite{PAR}).
So this critical point theorem still hold under the Cerami condition.
\begin{thm}\label{t2.1}(\cite{DL07})
Let
$\mathcal{H}$
be a real Banach space and let
$C_{-}$,
$C_{+}$
be two
symmetric cones in
$\mathcal{H}$
such that
$C_{+}$
is closed in
$\mathcal{H}$,
$C_{-}\cap C_{+} = \{0\}$
and
$$
i(C_{-}\setminus \{0\}) = i(\mathcal{H}\setminus C_{+}) = m < \infty.
$$
Define the following four sets by
\begin{eqnarray*}
&&D_{-}=\{u\in C_{-}|\,\|u\| \leq r_{-}\}, \\
&&S_{+}=\{u\in C_{+}|\,\|u\|=r_{+}\}, \\
&&Q=\{u + te|\, u\in C_{-},  t \geq 0, \|u +te\| \leq r_{-}\}, \;\; e\in \mathcal{H}\setminus C_-, \\
&&H=\{u + te|\,u\in C_{-},  t \geq 0,  \|u + te\| = r_{-}\}.
\end{eqnarray*}
Then
$(Q,  D_{-} \cup H)$
links
$S_{+}$
cohomologically in dimension
$m+1$
over
$\mathbf{Z}_{2}$.
Moreover,
suppose
$\Psi\in C^1(\mathcal{H},\mathbf{R})$
satisfying the Cerami condition,
and
$\displaystyle\sup_{x\in D_{-} \cup H}\Psi(x) < \displaystyle\inf_{x\in S^{+}}\Psi(x)$,
$\displaystyle\sup_{x\in Q}\Psi(x) < \infty$.
Then
$\Psi$
has a critical value
$d \geq \displaystyle\inf_{x\in S^{+}}\Psi(x)$.
\end{thm}

For convenience,
let us recall the definition and some properties
of the cohomological index of Fadell-Rabinowitz for a
$\mathbf{Z}_{2}$-set,
see \cite{FR77, FR78, PAR} for details.
For simplicity,
we only consider the usual
$\mathbf{Z}_{2}$-action on a
linear space,
i.e.,
$\mathbf{Z}_{2}=\{1, -1\}$
and the action is the
usual multiplication.
In this case,
the $\mathbf{Z}_{2}$-set
$A$
is a symmetric set with
$-A=A$.

Let
$E$
be a normed linear space.
We denote by
$\mathcal {S}(E)$
the set of all symmetric subsets of
$E$
which do not contain the
origin of
$E$.
For
$A\in \mathcal {S}(E)$,
denote
$\bar{A} = A/\mathbf{Z}_{2}$.
Let
$\rho: \bar{A} \rightarrow \mathbf{R}P^{\infty}$
be the classifying map and
$\rho^{*}: H^{*}(\mathbf{R}P^{\infty})=\mathbf{Z}_{2}[\omega] \rightarrow H^{*}(\bar{A})$
the induced homomorphism of the cohomology rings.
The cohomological index of $A$,
denoted by $i(A)$,
is defined by
$\displaystyle\sup\{k \geq 1: \rho^{*}(\omega^{k-1}) \neq 0\}$.
We list some properties of the cohomological index here for further use  in this paper.
Let
$A, B\in \mathcal {S}(E)$,
there hold
\begin{enumerate}
\item[(i1)]
({\bf monotonicity})
if
$h:A \rightarrow B$
is an odd map,
then $i(A)\le i(B)$,
\item[(i2)]
({\bf continuity})
if
$C$
is a closed symmetric subset of
$A$,
then there exists a closed  symmetric neighborhood
$N$
of
$C$
in
$A$,
such that
$i(N) = i(C)$,
hence the interior of
$N$
in
$A$
is also  a neighborhood of
$C$
in
$A$
and
$i({\rm int} N)=i(C)$,
\item[(i3)]
({\bf neighborhood of zero})
if
$V$
is bounded closed symmetric neighborhood of the origin in
$E$,
then
$i(\partial V) = \dim E$.
\end{enumerate}

\section{The eigenvalue problem}
First we solve the eigenvalue problem
\begin{equation}\label{e3.9}
E'({\bf u}) = \mu J'({\bf u}),\; {\bf u}\in \mathcal{H}.
\end{equation}
\begin{lem}\label{l3.1}
For any
${\bf u} = (u_{1}, u_{2}), {\bf v} = (v_{1}, v_{2})\in \mathcal{H}$,
it holds that
\begin{equation}\label{e3.17}
\langle E'({\bf u})-E'({\bf v}), {\bf u} - {\bf v}\rangle
\geq (\|u_{1}\|_{1}-\|v_{1}\|_{1})^2 + (\|u_{2}\|_{2}-\|v_{2}\|_{2})^2.
\end{equation}
\end{lem}
\pf.
By direct computations,
we have
\begin{eqnarray*}\begin{array}{lllllll}
\langle E'({\bf u})-E'({\bf v}), {\bf u} - {\bf v}\rangle \\
=\int_{\mathbf{R}^{N}}\left(|\nabla u_{1}|^{2}+|\nabla v_{1}|^{2}-2\nabla
u_{1}\cdot \nabla v_{1}\right){\rm d}x
+\int_{\mathbf{R}^{N}}b_{1}(x)
\left(|u_{1}|^{2}+|v_{1}|^{2}- 2u_{1}v_{1}\right){\rm d}x \\+ \int_{\mathbf{R}^{N}}\left(|\nabla u_{2}|^{2}+|\nabla v_{2}|^{2}-2\nabla
u_{2}\cdot \nabla v_{2}\right){\rm d}x
+\int_{\mathbf{R}^{N}}b_{2}(x)\left
(|u_{2}|^{2}+|v_{2}|^{2}- 2u_{2}v_{2}\right){\rm d}x.
\end{array}
\end{eqnarray*}
From the definition of the norm in
$H_{i}$,
we can get
\begin{equation}\label{e3.11}\begin{array}{ll}
\int_{\mathbf{R}^{N}}\left(|\nabla u_{1}|^{2}+|\nabla v_{1}|^{2}-2\nabla
u_{1}\cdot \nabla v_{1}\right){\rm d}x
+\int_{\mathbf{R}^{N}}b_{1}(x)
\left(|u_{1}|^{2}+|v_{1}|^{2}- 2u_{1}v_{1}\right){\rm d}x\\
= \|u_{1}\|_{1}^{2}+\|v_{1}\|_{1}^{2}
-2\langle u_1,\,v_1\rangle_1\geq
\|u_{1}\|_{1}^{2} + \|v_{1}\|_{1}^{2} - 2\|u_{1}\|_{1}\|v_{1}\|_{1}= (\|u_{1}\|_{1} - \|v_{1}\|_{1})^{2},\end{array}
\end{equation}
\begin{equation}\label{e3.16}\begin{array}{ll}
\int_{\mathbf{R}^{N}}\left(|\nabla u_{2}|^{2}+|\nabla v_{2}|^{2}-2\nabla
u_{2}\cdot \nabla v_{2}\right){\rm d}x
+\int_{\mathbf{R}^{N}}b_{2}(x)\left
(|u_{2}|^{2}+|v_{2}|^{2}- 2u_{2}v_{2}\right){\rm d}x\\
= \|u_{2}\|_{2}^{2}+\|v_{2}\|_{2}^{2}
-2\langle u_2,\,v_2\rangle_2\geq
\|u_{2}\|_{2}^{2} + \|v_{2}\|_{2}^{2} - 2\|u_{2}\|_{2}\|v_{2}\|_{2}= (\|u_{2}\|_{2} - \|v_{2}\|_{2})^{2}.\end{array}
\end{equation}
Now
(\ref{e3.11})
and
(\ref{e3.16})
imply
(\ref{e3.17}).
\qed
\begin{lem}\label{l3.2}
If
 ${\bf u}_{n}\rightharpoonup {\bf u}$
and
$\langle E'({\bf u}_{n}), {\bf u}_{n} - {\bf u}\rangle \rightarrow 0$,
then
${\bf u}_{n}\rightarrow {\bf u}$
in
$\mathcal{H}$.
\end{lem}
\pf.
Since
$\mathcal{H}$
is a Hilbert space and ${\bf u}_{n} = (u_{n}, v_{n})\rightharpoonup {\bf u} = (u, v)$,
we only need to show that
$\|{\bf u}_{n}\| \rightarrow \|{\bf u}\|$.
Note that
\begin{eqnarray*}
\lim_{n\rightarrow \infty}\langle E'({\bf u}_{n})-E'({\bf u}), {\bf u}_{n} - {\bf u}\rangle
=\lim_{n\rightarrow \infty}(\langle E'({\bf u}_{n}), {\bf u}_{n} - {\bf u}\rangle - \langle
E'({\bf u}), {\bf u}_{n} - {\bf u}\rangle)
= 0.
\end{eqnarray*}
By inequality
(\ref{e3.17})
we have
\begin{eqnarray*}
\langle E'({\bf u}_{n})-E'({\bf u}), {\bf u}_{n} - {\bf u}\rangle
\geq (\|u_{n}\|_{1}-\|u\|_{1})^2 + (\|v_{n}\|_{2}-\|v\|_{2})^2.
\end{eqnarray*}
So
$\|u_{n}\|_{1} \to \|u\|_{1}$,
$\|v_{n}\|_{2} \to \|v\|_{2}$
and hence
$\|{\bf u}_{n}\| \rightarrow \|{\bf u}\|$
as
$n \rightarrow \infty$
and the assertion follows.\qed
\begin{lem}\label{l3.3}
$J'$
is weak-to-strong continuous,
i.e.
${\bf u}_{n}\rightharpoonup {\bf u}$
in
$\mathcal{H}$
implies
$J'({\bf u}_{n})\rightarrow J'({\bf u})$.
\end{lem}
\pf.
Since
${\bf u}_{n} = (u_n, v_n)\rightharpoonup {\bf u} = (u,v)$
in
$\mathcal{H}$,
$u_n \rightharpoonup u$
in
$H_{1}$.
So
$u_n\to u$
in
$L^{2}(\mathbf{R}^{N})$
because
$H_{1}$
compactly embedded into
$L^{2}(\mathbf{R}^{N})$.
Similarly,
we have
$v_n\to v$
in
$L^{2}(\mathbf{R}^{N})$.
For any
${\bf v}=(\tilde u, \tilde v)\in \mathcal{H}$,
$$
\int_{\mathbf{R}^{N}}\tilde u^2{\rm d}x
\leq \frac{1}{b_{1}^{0}}\int_{\mathbf{R}^{N}}b_{1}(x)\tilde u^{2}{\rm d}x
\leq \frac{1}{b_{1}^{0}}\|\tilde u\|_{1}^{2}
\leq \frac{1}{b_{1}^{0}}\|{\bf v}\|^{2},
$$
so
$\left(\int_{\mathbf{R}^{N}}\tilde u^2{\rm d}x\right)^{\frac{1}{2}}\leq C\|{\bf v}\|$.
Similarly,
we have
$\left(\int_{\mathbf{R}^{N}}\tilde v^2{\rm d}x\right)^{\frac{1}{2}}\leq C\|{\bf v}\|$.
Then,
\begin{eqnarray*}
&&|\langle J'({\bf u}_{n}) - J'({\bf u}), {\bf v}\rangle|\\
&=&\left|\int_{\mathbf{R}^{N}}\left(V_{1}(x)(u_n - u)\tilde u + \gamma(x) (v_n -v)\tilde u + \gamma(x) (u_n - u)\tilde v+ V_{2}(x)(v_n -v)\tilde v\right){\rm d}x\right|\\
&\leq & \|V_{1}\|_{\infty}\left(\int_{\mathbf{R}^{N}}(u_n - u)^{2}{\rm d}x\right)^{\frac{1}{2}}\left(\int_{\mathbf{R}^{N}}\tilde u^{2}{\rm d}x\right)^{\frac{1}{2}} + \|\gamma\|_{\infty}\left(\int_{\mathbf{R}^{N}}(v_n -v)^{2}{\rm d}x\right)^{\frac{1}{2}}\left(\int_{\mathbf{R}^{N}}\tilde u^{2}{\rm d}x\right)^{\frac{1}{2}}\\
&&+ \|\gamma\|_{\infty}\left(\int_{\mathbf{R}^{N}}(u_n - u)^{2}{\rm d}x\right)^{\frac{1}{2}}\left(\int_{\mathbf{R}^{N}}\tilde v^{2}{\rm d}x\right)^{\frac{1}{2}} + \|V_{2}\|_{\infty}
\left(\int_{\mathbf{R}^{N}}(v_n -v)^{2}{\rm d}x\right)^{\frac{1}{2}}\left(\int_{\mathbf{R}^{N}}\tilde v^{2}{\rm d}x\right)^{\frac{1}{2}}\\
&\leq & C \left(\int_{\mathbf{R}^{N}}(u_n -u)^{2}{\rm d}x\right)^{\frac{1}{2}}\|{\bf v}\| +  C \left(\int_{\mathbf{R}^{N}}(v_n -v)^{2}{\rm d}x\right)^{\frac{1}{2}}\|{\bf v}\|
\to 0,
\end{eqnarray*}
hence
$J'({\bf u}_{n})\rightarrow J'({\bf u})$.
\qed
\begin{lem}\label{l3.4}
If
${\bf u}_{n}\rightharpoonup {\bf u}$
in
$\mathcal{H}$,
then
$J({\bf u}_{n}) \rightarrow J({\bf u})$.
\end{lem}
\pf.
\begin{eqnarray*}
2|J({\bf u}_{n}) - J({\bf u})|
&=& |\langle J'({\bf u}_{n}), {\bf u}_{n}\rangle - \langle J'({\bf u}), {\bf u}\rangle|\\
&=& |\langle J'({\bf u}_{n}) - J'({\bf u}),  {\bf u}_{n}\rangle + \langle J'({\bf u}), {\bf u}_{n} - {\bf u}\rangle|\\
&\leq & \|J'({\bf u}_{n}) - J'({\bf u})\| \|{\bf u}_{n}\| + o(1).
\end{eqnarray*}
Because
${\bf u}_{n} \rightharpoonup {\bf u}$,
${\bf u}_{n}$
is bounded.
From Lemma
\ref{l3.3},
we have
$J({\bf u}_{n}) \rightarrow J({\bf u})$.
\qed

In this section, we assume that $V_{1}$ and $V_{2}$ satisfy the following condition
$$
(**)\;\;\;\;meas\{x\in \mathbf{R}^{N}|\,V_{1}(x) > 0\} > 0\;\;
{\rm or}\;\;
meas\{x\in \mathbf{R}^{N}|\,V_{2}(x) > 0\} > 0.
$$
Set
$\mathcal{M} = \{{\bf u}\in \mathcal{H}|\, J({\bf u}) = 1\}$,
by
$(**)$,
we can see that
$\mathcal{M}$
is not empty,
see also Lemma
\ref{l3.7}
below.
Clearly,
$J({\bf u})=\frac{1}{2}\langle J'({\bf u}), {\bf u}\rangle$,
so
$1$
is a regular value of the functional
$J$.
Hence by the implicit theorem,
$\mathcal{M}$
is a
$C^{1}$-Finsler manifold.
It is complete,
symmetric,
since
$J$
is continuous and even.
Moreover,
$0$
is not contained in
$\mathcal{M}$,
so the trivial
$\mathbf{Z}_{2}$-action on
$\mathcal{M}$
is free.
Set
$\widetilde{E} = E|_{\mathcal{M}}$.
\begin{lem}\label{l3.5}
If
${\bf u}\in \mathcal{M}$
satisfies
$\widetilde{E}({\bf u}) = \mu$
and
$\widetilde{E}'({\bf u}) = 0$,
then
$(\mu, {\bf u})$
is a solution of
the functional equation (\ref{e3.9}).
\end{lem}
\pf.
By Proposition 3.54 in
\cite{PAR},
the norm of
$\widetilde{E}'({\bf u})\in T^{*}_{{\bf u}}\mathcal{M}$
is given by
$\|\widetilde{E}'({\bf u})\|_{{\bf u}}^{*}
= \displaystyle\min_{\nu\in \mathbf{R}}\|E'({\bf u}) - \nu J'({\bf u})\|^{*}$
(
here the norm
$\|\cdot\|^{*}_{{\bf u}}$
is the norm in the
fibre
$T^{*}_{{\bf u}}\mathcal{M}$,
and
$\|\cdot\|^{*}$
is the operator
norm,
the minimal can be attained was proved in Lemma 3.55 in
\cite{PAR}
).
Hence there exists
$\nu\in \mathbf{R}$
such that
$E'({\bf u}) - \nu J'({\bf u}) = 0$,
that is
$(\nu,  {\bf u})$
is a solution of the equation
(\ref{e3.9})
and
$
\mu
= \widetilde{E}({\bf u})
= \frac{1}{2}\langle E'({\bf u}), {\bf u}\rangle
= \frac{1}{2}\langle \nu J'({\bf u}), {\bf u}\rangle
= \frac{\nu}{2}\langle J'({\bf u}), {\bf u}\rangle
= \nu J({\bf u})
= \nu
$.
\qed
\begin{lem}\label{l3.6}
$\widetilde{E}$
satisfies the
$(PS)$
condition,
i.e.
if
$({\bf u}_{k})$
is a sequence on
$\mathcal{M}$
such that
$\widetilde{E}({\bf u}_{k}) \rightarrow c$,
and
$\widetilde{E}'({\bf u}_{k}) \rightarrow 0$,
then up to a subsequence
${\bf u}_{k} \rightarrow {\bf u}\in \mathcal{M}$
in
$\mathcal{H}$.
\end{lem}
\pf.
First,
from the definition of
$E$,
we can deduce that
$({\bf u}_{k})$
is bounded.
Then,
up to a subsequence,
${\bf u}_{k}$
converges weakly to some
${\bf u}$,
by Lemma \ref{l3.4},
we have
$J({\bf u})=1$,
so
${\bf u}\in \mathcal{M}$.

From
$\widetilde{E}'({\bf u}_{k}) \rightarrow 0$,
we have
$E'({\bf u}_{k}) - \nu_k J'({\bf u}_{k}) \rightarrow 0$
in
$\mathcal{H}$
for a sequence of real numbers
$(\nu_k)$.
So
$\langle E'({\bf u}_{k}) - \nu_k J'({\bf u}_{k}), {\bf u}_{k}\rangle \rightarrow 0$,
thus we get
$\nu_k \rightarrow c$.
By Lemma \ref{l3.3},
we have
$E'({\bf u}_{k})\rightarrow c J'({\bf u})$.
Hence
$$
\langle E'({\bf u}_{k}),  {\bf u}_{k} - {\bf u} \rangle
= \langle E'({\bf u}_{k})-cJ'({\bf u}),  {\bf u}_{k} - {\bf u} \rangle+\langle cJ'({\bf u}),  {\bf u}_{k} - {\bf u} \rangle \rightarrow 0.
$$
By Lemma \ref{l3.2},
we obtain
${\bf u}_{k} \rightarrow {\bf u}$.
\qed

Let
$\mathcal{F}$
denote the class of symmetric subsets of
$\mathcal{M}$,
$\mathcal{F}_{n}
= \{M\in \mathcal{F}|\, i(M)
\geq n\}$
and
\begin{equation}\label{e3.14}
\mu_{n} = \displaystyle\inf_{M\in
\mathcal{F}_{n}}\displaystyle\sup_{{\bf u}\in
M}{E}({\bf u}).
\end{equation}
Since
$\mathcal{F}_{n}
\displaystyle\supset \mathcal{F}_{n+1}$,
$\mu_{n}\leq \mu_{n+1}$.
\begin{lem}\label{l3.7}
If
$(**)$
holds,
then for every
$\mathcal{F}_{n}$,
there is a compact symmetric set
$M\in \mathcal{F}_{n}$.
\end{lem}
\pf.
We follow the idea of the proof of Theorem 3.2 in \cite{HT}.
Suppose
$meas\{x\in \mathbf{R}^{N}|\,V_{1}(x) > 0\} > 0$,
it implies  that
$\forall\, n\in \Bbb{N}$,
there exist
$n$
open balls
$(B_{i})_{1 \leq i \leq n}$
in
$\mathbf{R}^{N}$
such that
$B_{i}\cap B_{j} = \emptyset$
for
$i \neq j$
and
$meas(\{x\in\mathbf{R}^{N}|\, V_{1}(x) > 0\}\cap B_{i}) > 0$.
Approximating the characteristic function
$\chi_i$
of set
${\{x\in \mathbf{R}^{N}|\, V_{1}(x) > 0\}\cap B_{i}}$
by a
$C^{\infty}$-function
$u_i$
in
$L^{2}(\mathbf{R}^{N})$,
and require that the sequence
$\{u_i\}_{1 \leq i \leq n}\subseteq C^{\infty}(\mathbf{R}^{N})$
satisfies
$\int_{\mathbf{R}^{N}}V_{1}(x)|u_{i}|^{2}{\rm d}x > 0 $
for all
$i = 1,\cdots,n$
and
${\rm supp}\,u_{i} \cap  {\rm supp}\,u_{j}=\emptyset$
when
$i \neq j$.
Set
${\bf u}_{i} = (u_{i},0)\in \mathcal{H}$,
then
$J({\bf u}_{i})=\frac 12\int_{\mathbf{R}^{N}}V_{1}(x)|u_{i}|^{2}{\rm d}x > 0$.
Normalizing
${\bf u}_{i}$,
we assume that
$J({\bf u}_{i}) = 1$.
Denote by
$U_{n}$
the space spanned by
$({\bf u}_{i})_{1 \leq i \leq n}$.
$\forall {\bf u}\in U_{n}$,
we have
${\bf u} = \sum\limits_{i=1}^{n}\alpha_{i}{\bf u}_{i}$
and
$J({\bf u})=\sum\limits_{i=1}^{n}|\alpha_{i}|^{2}$.
So
$\left(J({\bf u})\right)^{\frac{1}{2}}$
defines a norm on
$U_{n}$.
Since
$U_{n}$
is
$n$
dimensional,
this norm is equivalent to
$\|\cdot\|$.
Thus
$\{{\bf u}\in U_{n}|\, J({\bf u})=1\}\subseteq \mathcal{M}$
is compact with respect to the norm
$\|\cdot\|$
and by the property (i3) of cohomological index,
$i(\{{\bf u}\in U_{n}|\,J({\bf u})=1\})=n$.
So
$\{{\bf u}\in U_{n}|\, J({\bf u})=1\}\in \mathcal{F}_{n}$.
If
$meas\{x\in \mathbf{R}^{N}|\,V_{2}(x) > 0\} > 0$,
the proof is similar.\qed

By Lemma \ref{l3.7},
we have $\mu_{n}<+\infty$,
and by condition (B),
there holds
$\mu_n\ge 0$.
Furthermore,
by Lemma \ref{l3.6} and Proposition 3.52 in \cite{PAR},
we see that
$\mu_n$
is sequence of critical values of
$\widetilde{E}$
and
$\mu_{n}\to +\infty,\; \mbox{as}\;n\to \infty$ .
By Lemma \ref{l3.5} we get a divergent sequence of eigenvalues for problem (\ref{e3.9}).
So we have the following result.
\begin{thm}\label{t3.2}
Under the condition $(**)$,
the problem (\ref{e3.9}) has an increasing sequence eigenvalues
$\mu_n$
which are defined by (\ref{e3.14}) and
$\mu_{n}\to +\infty,\; \mbox{as}\;n\to \infty$ .
\end{thm}
\begin{lem}\label{l3.8}
Under the condition
$(**)$,
Set
\begin{equation}\label{e3.15}
\rho_{n}
= \displaystyle\inf_{K\in \mathcal{F}^{c}_{n}}\displaystyle\sup_{{\bf u}\in K}E({\bf u}),
\end{equation}
where
$\mathcal{F}^{c}_{n}=\{K\in \mathcal{F}_n|\, K ~is~ compact \}$.
we have
$\mu_n=\rho_n$.
\end{lem}
\pf.
From Lemma \ref{l3.7},
$\mathcal{F}^{c}_{n}\ne \emptyset$
and so
$\rho_n<+\infty$.
It is obvious that
$\mu_n\le \rho_n$.
If
$\mu_n<\rho_n$,
there is
$M\in \mathcal{F}_{n}$
such that
$\sup\limits_{{\bf u}\in M}E({\bf u})<\rho_n$.
The closure
$\overline{M}$
of
$M$
in
$\mathcal{M}$
is still in
$\mathcal{F}_{n}$,
by continuity of
$E$,
$\sup\limits_{{\bf u}\in \overline{M}}E({\bf u})<\rho_n$
holds.
By the property (i2) of the cohomological index,
we can find a small open neighborhood
$A\in \mathcal{F}_{n}$ of $\overline{M}$
in
$\mathcal{M}$
such that
$\sup\limits_{{\bf u}\in A}E({\bf u})<\rho_n$.
As it was proved in the proof of Proposition
$3.1$
in
\cite{DL07},
for every symmetric open subset
$A$
of
$\mathcal{M}$,
there holds
$
i(A)
= \displaystyle\sup\{i(K)|\, K ~{\rm is ~compact~ and ~ symmetric~\\with} ~ K \subseteq A\}
$.
So we can choose a symmetric compact subset
$K\subseteq A$
with
$i(K)\ge n$
and
$\sup\limits_{{\bf u}\in K}E({\bf u})<\rho_n$.
This contradicts to the definition of
$\rho_n$.
Therefore we have
$\mu_n=\rho_n$.
\qed

Set
$C_m=\{{\bf u}\in \mathcal{H}\setminus \{0\}|\,E({\bf u}) \leq \mu_{m} J({\bf u})\}$
and
$D_m=\{{\bf u}\in \mathcal{H}|\, E({\bf u}) < \mu_{m+1} J({\bf u})\}$.
It is clear that
$C_m$,
$D_m\in \mathcal{S}(\mathcal{H})$,
i.e.,
$C_m$
and
$D_m$
are symmetric subsets of
$\mathcal{H}$
and do not contain 0.

\begin{thm}\label{t3.3}
If
$\mu_{m} < \mu_{m+1}$
for some
$m\in \Bbb N$,
then the cohomological indices satisfy
\begin{equation}\label{e3.aa}
i(C_m) = i(D_m)= m.
\end{equation}
\end{thm}
\pf.
Follow the idea of the proof of Theorem 3.2 in \cite{DL07}.
Suppose
$\mu_{m} < \mu_{m+1}$.
If we set
$A_m = \{{\bf u}\in \mathcal{M}|\, E({\bf u}) \leq \mu_{m}\}$
and
$B_m = \{{\bf u}\in \mathcal{M}|\,E({\bf u}) < \mu_{m+1}\}$,
by the definition (\ref{e3.14}),
we have
$i(A_m) \leq m$.
Assume that
$i(A_m) \leq m-1$.
Then,
by the property (i2) of the cohomological index,
there exists a symmetric neighborhood
$N$
of
$A_m$
in
$\mathcal{M}$
satisfying
$i(N) = i(A_m)$.
By the equivariant deformation theorem (see \cite{D03}),
there exists
$\delta > 0$
and an odd continuous map
$
\iota:\{{\bf u}\in \mathcal{M}|\,
E({\bf u}) \leq \mu_{m}+\delta\} \rightarrow
\{{\bf u}\in \mathcal{M}|\, E({\bf u}) \leq \mu_{m}-\delta\}\cup N = N
$.
Hence
$
i({\bf u}\in \mathcal{M}|\,E({\bf u})
\leq
\mu_{m}+\delta)
\leq m-1
$.
By (\ref{e3.14}),
there exists
$M\in \mathcal F_m$
such that
$\sup\limits_{{\bf u}\in M}E({\bf u})<\mu_m+\delta$.
So
$
M\subseteq \{{\bf u}\in \mathcal{M}|\,E({\bf u}) \leq
\mu_{m}+\delta\}
$
and thus
$i(M)\le m-1$.
This contradicts to the fact that
$M\in \mathcal F_m$.
Thus we have
$i(A_m) = m$.
By  $2$-homogeneousness of the functionals
$E, J$,
the map
$h: C_m\to A_m$
with
$h({\bf u})=\frac{{1}}{\sqrt{J({\bf u})}}{\bf u}$
is odd,
from the monotonicity (i1) of the cohomological index,
we have
$i(C_m) \le m$.
But it is clear that
$A_m\subset C_m$,
we have
$i(C_m) \ge m$,
so
$i(C_m) = m$.

Since
$A_m\subseteq B_m$
and
$i(A_m)=m$,
we have
$i(B_m) \geq m$.
Assume that
$i(B_m) \geq m+1$.
As in the proof of Lemma \ref{l3.8},
there exists a symmetric,
compact subset
$K$
of
$B_m$
with
$i(K) \geq m+1$.
Since
$\max\limits_{ {\bf u}\in K}E({\bf u})< \mu_{m+1}$,
this contradicts to definition (\ref{e3.14}).
So
$i(B_m)=m$.
Similar to the above arguments,
we also have
$i(D_m) = m $.
\qed
\begin{rem}\label{re3.1}
{\rm
If we consider the following eigenvalue problem,
\begin{equation}\label{e3.18}
E'({\bf u}) = \mu J'({\bf u}),\; {\bf u}\in \mathcal{H}_{r},
\end{equation}
then  all the results in this section still hold,
 we only need to replace the space
$\mathcal{H}$
  by
$\mathcal{H}_{r}$.

}
\end{rem}

\section{Proof of the main theorems}
Replacing
$(\lambda,  V_i, \gamma)$
with
$(-\lambda,  -V_i, -\gamma)$ if necessary,
we can assume that
$\lambda \geq 0$.
First,
we consider the case that condition
($**$) holds
and there exists
$m \geq 1$
such that
$\mu_{m} \leq \lambda <\mu_{m+1}$.
Set
\begin{equation}\label{e4.16}
C_{-}
= \{{\bf u}\in \mathcal{H}|\, E({\bf u})
\leq \mu_{m}J({\bf u})\},
\end{equation}
\begin{equation}\label{e4.17}
C_{+}
= \{{\bf u}\in \mathcal{H}|\, E({\bf u})
\geq \mu_{m+1} J({\bf u})\}.
\end{equation}
It is easy to see that
$C_{-}$,
$C_{+}$
are two symmetric closed cones in
$\mathcal{H}$
and
$C_{-} \cap C_{+} = \{0\}$.
By (\ref{e3.aa}) we have
\begin{equation}\label{e4.18}
i(C_{-}\setminus\{0\})
= i(C_m)
= i(D_m)
= i(\mathcal{H}
\setminus C_{+})
= m.
\end{equation}
\begin{lem}\label{t4.4}
There exist
$r_{+} > 0$
and
$\alpha > 0$
such that
$\Psi({\bf u}) > \alpha$
for
${\bf u}\in C_{+}$
and
$\|{\bf u}\| = r_{+}$.
\end{lem}
\pf.
Let
$\varepsilon > 0$
be small enough,
from (W$_1$) and (W$_3$),
we have
$|W(x, z)| \leq \varepsilon|z|^{2}+C_{\varepsilon}|z|^{p}$.

By the Sobolev embedding inequality,
for
${\bf u} = (u_{1}, u_{2})\in C_+$,
we can get
\begin{equation}\label{e4.18}
\begin{array}{lllllllll}
$$
\Psi({\bf u})
&=& E({\bf u})-\lambda J({\bf u})- P({\bf u}) \\
&=& E({\bf u}) - \frac{\lambda}{\mu_{m+1}}\cdot\mu_{m+1}J({\bf u})-P({\bf u})\\
&\geq& E({\bf u}) - \frac{\lambda}{\mu_{m+1}}E({\bf u}) -
\varepsilon \int_{\mathbf{R}^{N}}|u_{1}|^{2}{\rm d}x \\&&- \varepsilon \int_{\mathbf{R}^{N}}|u_{2}|^{2}{\rm d}x -
C_{\varepsilon}\int_{\mathbf{R}^{N}}|u_{1}|^{p}{\rm d}x - C_{\varepsilon}\int_{\mathbf{R}^{N}}|u_{2}|^{p}{\rm d}x \\
&\geq& E({\bf u}) - \frac{\lambda}{\mu_{m+1}}E({\bf u}) -
\frac{\varepsilon}{b_{1}^{0}} \int_{\mathbf{R}^{N}}b_{1}(x)|u_{1}|^{2}{\rm d}x - \frac{\varepsilon}{b_{2}^{0}} \int_{\mathbf{R}^{N}}b_{2}(x)|u_{2}|^{2}{\rm d}x \\&&-
C_{\varepsilon}\int_{\mathbf{R}^{N}}|u_{1}|^{p}{\rm d}x - C_{\varepsilon}\int_{\mathbf{R}^{N}}|u_{2}|^{p}{\rm d}x\\
&\geq& (1 - \frac{\lambda}{\mu_{m+1}} -  2\max(\frac{\varepsilon}{b_{1}^{0}}, \frac{\varepsilon}{b_{2}^{0}}))E({\bf u})- C_{\varepsilon}\int_{\mathbf{R}^{N}}|u_{1}|^{p}{\rm d}x - C_{\varepsilon}\int_{\mathbf{R}^{N}}|u_{2}|^{p}{\rm d}x\\
&\geq& \frac{1}{2}(1-\frac{\lambda}{\mu_{m+1}}-2\max(\frac{\varepsilon}{b_{1}^{0}},\frac{\varepsilon}{b_{2}^{0}}))\|{\bf u}\|^{2} - C\|{\bf u}\|^{p}.
\end{array}
\end{equation}
We remind that in the second inequality of (\ref{e4.18}),
the condition (B) has been applied.
Since
$p>2$,
the assertion follows.\qed

Since
$\lambda \geq \mu_{m}$,
by (W$_1$) it holds that
\begin{equation}\label{e4.19}
\Psi({\bf u}) \leq 0,\;\forall\, {\bf u}\in C_{-}.
\end{equation}
Set
$\mathbf R^+=[0,+\infty)$.
Following the idea of the proof of Theorem 4.1 in \cite{DL07},
we have
\begin{lem}\label{t4.5}
Let
${\bf e} = (e_{1}, e_{2})\in \mathcal{H}\setminus C_{-}$,
there exists
$r_{-} > r_{+}$
such that
$\Psi({\bf u}) \leq 0$
for
${\bf u}\in C_{-} + \mathbf{R}^{+}{\bf e}$
and
$\|{\bf u}\| \geq r_{-}$.
\end{lem}
\pf.
Define another norm on
$\mathcal{H}$
by
$
\|{\bf u}\|_{V}^2:=\int_{\mathbf{R}^{N}}(|V_{1}(x)|+|\gamma(x)| + 1)|u|^{2}{\rm d}x
+ \int_{\mathbf{R}^{N}}(|V_{2}(x)| + |\gamma(x)| + 1)|v|^{2}{\rm d}x
$
for
${\bf u}=(u,v)$.
Then the same reason as the proof of Theorem 4.1 in \cite{DL07},
there exists some constant
$b > 0$
such that
$\|{\bf u} + t{\bf e}\| \leq b \|{\bf u} + t{\bf e}\|_{V}$
for every
${\bf u}\in C_{-}$,
$t \geq 0$
and some
$b > 0$.
That is
\begin{equation}\label{e4.20}
\begin{array}{ll}
\int_{\mathbf{R}^{N}}(|\nabla (u + te_{1})|^{2} + b_{1}(x)|u+ te_{1}|^{2}){\rm d}x +
\int_{\mathbf{R}^{N}}(|\nabla (v + te_{2})|^{2} + b_{2}(x)|v+ te_{2}|^{2}){\rm d}x \\
\leq b^{2}\int_{\mathbf{R}^{N}}(|V_{1}(x)|+|\gamma(x)|+1)|u + te_{1}|^{2}{\rm d}x + b^{2}\int_{\mathbf{R}^{N}}(|V_{2}(x)|+|\gamma(x)|+1)|v + te_{2}|^{2}{\rm d}x.
\end{array}
\end{equation}

Let
$\{{\bf u}_{k}\}$
be a sequence such that
$\|{\bf u}_{k}\|\rightarrow +\infty$
and
${\bf u}_{k}\in C_{-}+\mathbf{R}^{+}{\bf e}$.
Set
${\bf v}_{k}=(u_k,v_k):=\frac{{\bf u}_{k}}{\|{\bf u}_{k}\|}$,
then,
up to a subsequence,
$\{{\bf v}_{k}\}$
converges to some
${\bf v}=(u_0,v_0)$
weakly in
$\mathcal{H}$
and
$u_k\to u_0$,
$v_k\to v_0$
a.e.
in
$\mathbf{R}^{N}$.
Note that Lemma \ref{l3.4} is also true for functional
$
\int_{\mathbf{R}^{N}}(|V_{1}(x)|+|\gamma(x)|+1)|u|^{2}{\rm d}x + \int_{\mathbf{R}^{N}}(|V_{2}(x)|+|\gamma(x)|+1)|v|^{2}{\rm d}x,\;{\bf u}=(u,v)\in \mathcal{H}
$,
it follows from (\ref{e4.20}) that
$
\int_{\mathbf{R}^{N}}(|V_{1}(x)|+|\gamma(x)|+1)|u_0|^{2}{\rm d}x +\int_{\mathbf{R}^{N}}(|V_{2}(x)|+|\gamma(x)|+1)|v_0|^{2}{\rm d}x \geq \frac{1}{b^{2}}
$.
So
$|{\bf v}| \neq 0$
on a  positive measure set
$\Omega_0$.
By (W$_2$) we have
$$
\lim_{k \rightarrow \infty}\frac{W(x, {\bf u}_k(x))}{\|{\bf u}_{k}\|^{2}}
=\lim_{k\rightarrow\infty}\frac{W(x, \|{\bf u}_{k}\|{\bf v}_k(x))}
{\|{\bf u}_{k}\|^{2}|{\bf v}_{k}(x)|^{2}} |{\bf v}_{k}(x)|^{2}
=
+\infty,
\;x\in
\Omega_0.
$$
By (W$_1$) and Fatou lemma
we can get
$$
\frac{\int_{\mathbf{R}^{N}}W(x, {\bf u}_k(x))
{\rm d}x}{\|{\bf u}_{k}\|^{2}}\rightarrow +\infty,
\;{\rm as}\;k
\to
\infty.
$$
By the arbitrariness of the sequence
$\{{\bf u}_{k}\}$,
we have
\begin{equation}\label{e4.a}
\frac{\int_{\mathbf{R}^{N}}W(x, {\bf u}(x))
{\rm d}x}{\|{\bf u}\|^{2}}\rightarrow +\infty
\end{equation}
as
$\|{\bf u}\|\rightarrow +\infty$
and
${\bf u}\in C_{-}+\mathbf{R}^{+}{\bf e}$.
Noting that
\begin{equation}\label{e4.b}
\frac{\Psi({\bf u})}{\|{\bf u}\|^{2}} = \frac{1}{2} - \frac{\lambda J({\bf u})}{\|{\bf u}\|^{2}} -
\frac{\int_{\mathbf{R}^{N}}W(x, {\bf u}(x))
{\rm d}x}{\|{\bf u}\|^{2}}
\end{equation}
and by conditions (B) and (V), for
${\bf u}=(u,v)\in \mathcal{H}$
\begin{equation}\label{e4.c}
\left|\frac{J({\bf u})}{\|{\bf u}\|^{2}}\right| \leq \frac{C(\int_{\mathbf{R}^{N}}|u|^{2}{\rm d}x + \int_{\mathbf{R}^{N}}|v|^{2}{\rm d}x)}{\|{\bf u}\|^{2}}
\leq \frac{C(\int_{\mathbf{R}^{N}}b_{1}(x)|u|^{2}{\rm d}x +
\int_{\mathbf{R}^{N}}b_{2}(x)|v|^{2}{\rm d}x)}{\|{\bf u}\|^{2}}\leq C,
\end{equation}
the assertion follows from
(\ref{e4.a}),
(\ref{e4.b})
and
(\ref{e4.c}).
\qed

\begin{lem}\label{t4.6}
$\Psi$ satisfies the Cerami condition,
i.e.,
for any sequence
$\{{\bf u}_{k}\}$
in
$\mathcal{H}$
satisfying
$(1 + \|{\bf u}_{k}\|)\Psi'({\bf u}_{k}) \rightarrow 0$
and
$\Psi({\bf u}_{k})\rightarrow c$
possesses a convergent subsequence.
\end{lem}
\pf.
Let
$\{{\bf u}_{k}\}$
be a sequence in
$\mathcal{H}$
satisfying
$(1 + \|{\bf u}_{k}\|)\Psi'({\bf u}_{k}) \rightarrow 0$
and
$\Psi({\bf u}_{k})\rightarrow c$.
We claim that
$\{{\bf u}_{k}\}$
is bounded in
$\mathcal{H}$.
Otherwise,
if
$\|{\bf u}_{k}\|\to \infty$,
we consider
${\bf v}_{k}:=\frac{{\bf u}_{k}}{\|{\bf u}_{k}\|}$.
Then,
up to subsequence,
we get
${\bf v}_k \rightharpoonup {\bf v}$ in $\mathcal{H}$
and
${\bf v}_k \to {\bf v}$
a.e.
in
$\mathbf{R}^{N}$.

If
$ {\bf v} \neq 0$
in
$\mathcal{H}$,
since
$\Psi'({\bf u}_{k}){\bf u}_{k} \rightarrow 0$,
that is to say
\begin{equation}\label{e4.dd}
\begin{array}{ll}
\|{\bf u}_k\|^2- \lambda J'({\bf u}_k)\cdot {\bf u}_k- \int_{\mathbf{R}^{N}}\nabla_z W(x, {\bf u}_k(x))\cdot {\bf u}_k{\rm d}x\\
=\|{\bf u}_k\|^2- 2\lambda J ({\bf u}_k)-
\int_{\mathbf{R}^{N}}\nabla_z W(x, {\bf u}_k(x))\cdot {\bf u}_k{\rm d}x \rightarrow 0,
\end{array}
\end{equation}
from
(\ref{e4.c}),
we have
$\frac{|J({\bf u}_k)|}{\|{\bf u}_k\|^2} \leq C$,
so by dividing the left hand side of
(\ref{e4.dd})
with
$\|{\bf u}_{k}\|^2$
there holds
\begin{equation}\label{e4.21}
\left|
\int_{\mathbf{R}^{N}}\frac{\nabla_z W(x, {\bf u}_k(x))\cdot {\bf u}_k(x)}{\|{\bf u}_{k}\|^{2}}{\rm d}x
\right|
\leq C'
\end{equation}
for some constant
$C'>0$.
On the other hand,
Since
${\bf v}(x)\ne 0$
in some positive measure set
$\Omega\subset \mathbf{R}^{N}$,
so
${\bf v}_k(x)\ne 0$
for large
$k$,
and
$|{\bf u}_k(x)|\to +\infty$
as
$k\to \infty$,
for any fixed
$x\in \Omega$.
So by
(W$_2$),
we have
\begin{equation}\label{e4.pp}
\lim_{k\to \infty}|{\bf v}_{k}(x)|^{2}\frac{2W(x,{\bf u}_k(x))}{|{\bf u}_{k}|^{2}}= +\infty,\;\; \forall\,x\in \Omega.
\end{equation}
By Remark (1) before Theorem \ref{t1.1},
we have
$$
\nabla_z W(x, {\bf u}_k(x))\cdot {\bf u}_k(x)\ge 2W(x,{\bf u}_k(x)).
$$
So as
$k\to +\infty$,
we have
$$
\int_{\mathbf{R}^{N}}\frac{\nabla_z W(x, {\bf u}_k(x))\cdot {\bf u}_k(x)}{\|{\bf u}_{k}\|^{2}}{\rm d}x
=\int_{\{{\bf v}_{k}(x) \neq 0\}}|{\bf v}_{k}(x)|^{2}\frac{\nabla_z W(x, {\bf u}_k(x))\cdot {\bf u}_k(x)}{|{\bf u}_{k}(x)|^{2}}{\rm d}x
$$
$$
\ge \int_{\mathbf{R}^{N}}\chi_{\{{\bf v}_{k}\ne 0\}}(x)|{\bf v}_{k}(x)|^{2}\frac{2W(x,{\bf u}_k(x))}{|{\bf u}_{k}(x)|^{2}}{\rm d}x
\ge
\int_{\Omega}\chi_{\{{\bf v}_{k}\ne 0\}}(x)|{\bf v}_{k}(x)|^{2}\frac{2W(x,{\bf u}_k(x))}{|{\bf u}_{k}(x)|^{2}}{\rm d}x
\to
\infty,
$$
this contradicts to (\ref{e4.21}).
There is another explanation about the above estimate.
We observe that there exists
$\delta>0$
such that
$meas\{x\in \Omega|\,|{\bf v}(x)|\ge \delta\}>0$.
Otherwise,
$\forall\,n\in \Bbb N$,
$meas\{x\in \Omega|\,|{\bf v}(x)|\ge \frac 1n \}=0$.
Set
$\Omega_n=\{x\in \Omega|\,|{\bf v}(x)|\ge \frac 1n \}$,
then in
$\Omega\setminus \bigcup_{n=1}^{+\infty}\Omega_n$,
there holds
${\bf v}(x)=0$.
But
$\Omega\setminus \bigcup_{n=1}^{+\infty}\Omega_n$
and
$\Omega$
have the same measure,
it is impossible.
We may assume
$meas\;\Omega<+\infty$,
by Egorov's theorem,
there exists a positive measure subset
$\Omega_0$
of
$\{x\in \Omega|\,|{\bf v}(x)|\ge \delta\}$
such that
${\bf v}_k$
uniformly convergent to
${\bf v}$,
so for
$k\ge K$
with
$K$
large, there holds
$|{\bf v}_k(x)|\ge \delta/2$
in
$\Omega_0$.
Thus
$(\ref{e4.pp})$
holds in
$\Omega_0$.
So there holds
$$
\int_{\{{\bf v}_{k}(x) \neq 0\}}|{\bf v}_{k}(x)|^{2}\frac{\nabla_z W(x, {\bf u}_k(x))\cdot {\bf u}_k(x)}{|{\bf u}_{k}(x)|^{2}}{\rm d}x \ge \int_{\Omega_0}|{\bf v}_{k}(x)|^{2}\frac{2W(x,{\bf u}_k(x))}{|{\bf u}_{k}(x)|^{2}}{\rm d}x
\to
\infty.
$$

If
${\bf v} = 0$
in
$\mathcal{H}$,
inspired by \cite{J99},
we choose
$t_{k}\in [0,1]$
such that
$\Psi(t_{k}{\bf u}_{k}):= \displaystyle\max_{t\in [0, 1]}\Psi(t{\bf u}_{k})$.
For any
$\beta> 0$
and
$\tilde{{\bf v}}_{k}:=(4\beta)^{1/2}{\bf v}_{k}\rightharpoonup 0$,
by Lemma \ref{l3.4} and the compactness of
$P'$
(see Lemma 1.22 in \cite{ZS06})
we have that
$J(\tilde{{\bf v}}_{k})\to 0$
and
$
\int_{\mathbf{R}^{N}}
W(x, \tilde{{\bf v}}_{k}(x))
{\rm d}x=P(\tilde{{\bf v}}_{k})-P(0)=\langle P'(\xi_k\tilde{{\bf v}}_{k}), \tilde{{\bf v}}_{k}\rangle=\langle P'(\xi_k\tilde{{\bf v}}_{k})-P'(0), \tilde{{\bf v}}_{k}\rangle+\langle P'(0), \tilde{{\bf v}}_{k}\rangle\to 0
$
as
$k\to \infty$,
here
$\xi_k\in (0,1)$.
So there holds
$$
\Psi(t_{k}{\bf u}_{k}) \geq \Psi(\tilde{{\bf v}}_{k}) =
2\beta - \lambda J(\tilde{{\bf v}}_{k}) - \int_{\mathbf{R}^{N}}
W(x, \tilde{{\bf v}}_{k}(x))
{\rm d}x
\geq \beta,
$$
when
$k$
is large enough.
By the arbitrariness of
$\beta$,
it implies that
\begin{equation}\label{e4.22}
\lim\limits_{k \rightarrow \infty}\Psi(t_{k}{\bf u}_{k}) = \infty.
\end{equation}
Since
$\Psi(0)=0,\;\Psi({\bf u}_{k})\to c$,
we have
$t_{k}\in (0, 1)$.
By the definition of
$t_{k}$,
\begin{equation}\label{e4.23}
\langle \Psi'(t_{k}{\bf u}_{k}), t_{k}{\bf u}_{k}\rangle = 0.
\end{equation}
From
(\ref{e4.22}),
(\ref{e4.23}),
we have
\begin{equation}\label{e4.24}
\begin{array}{lll}
\Psi(t_{k}{\bf u}_{k})-\frac{1}{2}\langle \Psi'(t_{k}{\bf u}_{k}), t_{k}{\bf u}_{k}\rangle \\
=\int_{\mathbf{R}^{N}}
\left(
\frac{1}{2}\nabla_{z}W(x, t_{k}{\bf u}_{k}(x))\cdot t_{k}{\bf u}_{k}(x)
- W(x, t_{k}{\bf u}_{k}(x))
\right)
{\rm d}x
\rightarrow
\infty.
\end{array}
\end{equation}
By
(W$_4$),
there exists
$\theta\ge 1$
such that
\begin{equation}\label{e4.25}
\begin{array}{ll}\int_{\mathbf{R}^{N}}
\left(
\frac{1}{2}\nabla_{z}W(x, {\bf u}_{k}(x))\cdot {\bf u}_{k}(x)
- W(x, {\bf u}_{k}(x))
\right)
{\rm d}x\\
\geq
\frac{1}{\theta}
\int_{\mathbf{R}^{N}}
\left(\nabla_{z}W(x, t_{k}{\bf u}_{k}(x))\cdot t_{k}{\bf u}_{k}(x)
- W(x, t_{k}{\bf u}_{k}(x))
\right)
{\rm d}x,\end{array}
\end{equation}
Hence
\begin{equation}\label{e4.27}
\int_{\mathbf{R}^{N}}
\left(
\frac{1}{2}\nabla_{z}W(x, {\bf u}_{k}(x))\cdot {\bf u}_{k}(x)
- W(x, {\bf u}_{k}(x))
\right)
{\rm d}x
\rightarrow
\infty.
\end{equation}
On the other hand,
\begin{equation}\label{e4.28}
\int_{\mathbf{R}^{N}}
\left(
\frac{1}{2}\nabla_{z}W(x, {\bf u}_{k}(x))\cdot {\bf u}_{k}(x)
- W(x, {\bf u}_{k}(x))
\right)
{\rm d}x
= \Psi({\bf u}_{k}) - \frac{1}{2}\langle\Psi'({\bf u}_{k}),{\bf u}_{k}\rangle
\to c.
\end{equation}
(\ref{e4.27})
and
(\ref{e4.28})
are contradiction.
Hence
$\{{\bf u}_{k}\}$
is bounded in
$\mathcal{H}$.
So up to a subsequence,
we can assume that
${\bf u}_{k}\rightharpoonup {\bf u}$
for some
$\mathcal{H}$.

Since
$\Psi'({\bf u}_{k}) = E'({\bf u}_{k})- \lambda J'({\bf u}_{k}) - P'({\bf u}_{k})\to 0$
and
$J'$,
$P'$
are compact,
we have that
$E'({\bf u}_{k})\to \lambda J'({\bf u})+ P'({\bf u})$
in
$\mathcal{H}$.
So
$$
\langle E'({\bf u}_{k}),{\bf u}_{k}-{\bf u}\rangle
=
\langle E'({\bf u}_{k})-(\lambda J'({\bf u})+ P'({\bf u})),{\bf u}_{k}-{\bf u}\rangle+\langle \lambda J'({\bf u})+ P'({\bf u}),{\bf u}_{k}-{\bf u}\rangle
\to
0.
$$
By Lemma \ref{l3.2},
${\bf u}_{k}\to u$
in
$\mathcal{H}$.
\qed

\begin{rem}\label{re4.1}
{\rm
If we replace the space
$\mathcal{H}$
by
$\mathcal{H}_{r}$,
then Lemma
\ref{t4.4},
\ref{t4.5},
\ref{t4.6}
also hold.
}
\end{rem}

\noindent{\bf Proof of Theorem \ref{t1.1}}
Define
$D_{-}$,
$S_{+}$,
$Q$,
$H$
as Theorem \ref{t2.1},
then from Lemma \ref{t4.4},
$\Psi({\bf u}) \geq \alpha > 0$
for every
${\bf u}\in S_{+}$,
from Lemma \ref{t4.5},
$\Psi({\bf u}) \leq 0$
for every
${\bf u}\in D_{-}\cup H$
and
$\Psi$
is bounded on
$Q$.
Applying Lemma \ref{t4.6},
it follows from Theorem \ref{t2.1} that
$\Psi$
has a critical value
$d \geq \alpha > 0$.
Hence
${\bf u}$
is a non-trivial weak solution of (\ref{e1.1}).

For the cases
$0 \leq \lambda < \mu_{1}$
or
$V_{1}^{+}(x) \equiv 0 \equiv V_{2}^{+}(x)$,
set
$C_{-} = \{0\}$
and
$C_{+} = \mathcal{H}$,
it is easy to see that the arguments above are also valid.
The proof of Theorem \ref{t1.1} is complete.
\qed

\noindent{\bf Proof of Theorem \ref{t1.2}}
By Remarks \ref{re3.1} and \ref{re4.1},
the proof is the same as that of Theorem \ref{t1.1},
we only need to replace the space
$\mathcal{H}$
 by
$\mathcal{H}_{r}$.
\qed

\end{document}